\documentclass[a4paper,12pt,notitlepage]{report}
\usepackage{amsmath,amsfonts,amssymb,color}
\begin{document}
\author{Amalaswintha Wolfsdorf}
\title{A.M.D.G.\\
Factorising Polynomials over Finite Fields}
\date{March 2006}
\maketitle
\begin{abstract}

The aim of this paper is to show that there exists a deterministic
algorithm that can be applied to compute the factors of a polynomial
of degree 2, defined over a finite field, given certain conditions.

\end{abstract}

\tableofcontents
\newpage
\thispagestyle{plain}
\section{Preface}

The study of prime numbers has been puzzling Number Theorists for
several centuries. Certainly since the remarkable results found by
Pierre de Fermat in the $17^{th}$ century, a new wave
 of motivation has triggered some of the most talented mathematicians
 to research this field in more depth.

The field is vast, and possibly one of the most challenging ones:
whilst the statement of a problem in this field may at first sound
like a lunchtime brainteaser for a hobby - number-cruncher, its
solution will in general be extremely complex and in many cases has
taken centuries
to find, if this has been achieved at all yet!\\

But this field is not only of high importance to theoretical
research. The most recent developments in this field have been
concentrated on computational number theory. The fact that still so
little is known about prime numbers, and that it is such a difficult
field to make much progress in, has been exploited by the computer
industry during the last century.

Secure transmission of data is made possible by prime numbers, and
hence research in Number Theory is nowadays mainly revolved around
finding ways to ensure that this level of security is maintained.\\

At the heart of this lies the problem of finding roots of
polynomials modulo prime numbers.

Whilst it is in theory \textit{possible} to do this, the procedures
that we know about so far are not very efficient and would in
general
take far too long to be of any practical use.\\

This dissertation (unfortunately) does not provide us with a magic
key to cracking such codes.  I will show that there exists a
deterministic algorithm that can, under certain circumstances, find
the factors of polynomials modulo a prime number. However the
running time of this algorithm is still much higher than some
probabilistic (and fairly reliable!) algorithms that are in use
already.\\
\nopagebreak
 The result that we will obtain here is hence
rather of interest to mathematicians working in Algorithmic Number
Theory than of practical use. Perhaps, however, similar techniques
will eventually be developed that might be put to more use in
practice. Perhaps the purely theoretical side of mathematics will
find its applications in practice, and Albert Einstein will be
proved wrong for his remark ``As far as the laws of mathematics
refer to reality, they are not certain; and as far as they are
certain, they do not refer to reality.''

\chapter{Introduction} \label{ch:intro}

The idea for this project originates from a claim made by Dr Neeraj
Kayal in
2005, together with some further refinements added by Prof Bjorn Poonen (University of California, Berkeley). \\

A well-known open problem in Algorithmic Number Theory is the
efficient calculation of roots of polynomials modulo a prime number
$q$ in deterministic polynomial-time.\\

A very basic example of this is the following: let $q$ be prime and
$a$ and number between $0$ and $q-1$. The study of Elementary Number
Theory provides us with easy tools to check whether there exists a
number $b$ between $0$ and $q-1$ such that $b^2 = a$ (mod $q$) - in
that case, $a$ is called a \textit{quadratic residue}.\\

For example, we could apply what is called ``Euler's Criterion''. It
says that if $q$ is an odd prime, then for all $a \in \mathbb{N}$,
we have

\begin{displaymath} \bigg(\frac{a}{q}\bigg) \equiv a^{(q-1)/2}
\textrm{ (mod } q) \, ,\end{displaymath}

where the fraction on the left hand side denotes the Legendre
Symbol, defined by

\begin{equation*}\bigg(\frac{a}{q}\bigg): = \left\{
\begin{array}{ll}
1 &
\textrm{if such a number $b$ exists} \\
-1 & \textrm{if no such $b$ exists} \\
0 & \textrm{if $q$ divides $a$}
\end{array} \right. \end{equation*}

But how can we calculate this number $b$, if it exists? We would
need to solve the equation $h(z):= z^2 - a \equiv 0$ modulo $q$.
This is a much harder problem, if it is to be solved efficiently.
Of course we could try substituting every value in
$\{0,1,\dots,q-1\}$ for $z$ to check whether the equation is
satisfied; however as we are in practical applications more
concerned
with large primes, this could take quite a while.\\

\section{The Claim}

Kayal claimed that we can factorise such a polynomial $h(z)$ defined
over a finite field $\mathbb{F}_q$ using a deterministic algorithm
with running time bounded by a universal polynomial in $(log \, q)$,
given certain circumstances (Poonen's input to this claim will be
discussed later). Roughly speaking, the underlying condition is that
we can construct an algebraic family of bivariate polynomials
$C_z(X, Y)$, each member of
which has a different number of solutions modulo $q$.\\

I will restrict the detailed proof to the case $deg \, (h) = 2$. The
case for higher degree polynomials will be discussed briefly
afterwards.\\
 I will then also give a brief discussion about the
running
time of the algorithm.\\

\section{The idea of the proof}

The idea of the proof can be outlined as follows:\\

We have a deterministic algorithm, known as ``Schoof's Algorithm'',
which is used to compute the number of rational points on
an elliptic curve given in Weierstrass Form and defined over a finite field.\\

Let $h(z)$ denote the polynomial that we wish to factorise, and
$\mathbb{F}_q$ the finite field
over which $h$ is defined.  We consider the ring $R:= \mathbb{F}_q[z]/(h(z))$.\\

If we are given an elliptic curve $C$ over $R$ that satisfies the
underlying condition, and we attempt to apply Schoof's Algorithm to
count the number of rational points on $C$, the algorithm will at
some point break down, and thereby reveal the factors of $h(z)$.\\

Now why does this happen?

Consider the difference between a ring and a field:  a field
contains all its inverses, which is not necessarily true for a ring.
This is precisely why the algorithm will not work when it is working
with a ring: whilst trying to compute the inverse of an element
(which the algorithm can easily do when in a field), it will at some
point not be able to find that inverse and will therefore stop running.\\

At this point we know that it must have found an element of the ring
that has no inverse.  But by inspecting this ring $R$ more closely,
we can see which elements in the ring do not have an inverse:  it is
precisely the set of elements in $\mathbb{F}_q[z]$
spanned by the factors of $h(z)$.\\

So all we need to do is compute the greatest common divisor of this
element that made the algorithm stop, and $h(z)$ (since this element
may be a multiple of a factor of a factor), to obtain a non-trivial
factor of $h(z)$!\\

To present a detailed proof however requires a lot more careful
explanation; this is what will follow now.\\

\section{Structure of this Paper}

In Chapter (\ref{ch:ecs}) of this dissertation, I will define
elliptic curves and explain some of their elementary properties that
we will need to be aware of in order to understand Schoof's
Algorithm.

Chapter (\ref{ch:counting}) contains a brief discussion about
counting rational points on elliptic curves, which is followed by a
rather technical section explaining the essential tools that
underlie
Schoof's Algorithm.\\

I will give a full description of Schoof's Algorithm for elliptic
curves over a finite field in Chapter (\ref{ch:schoof}).

For the purpose of a clear and thorough understanding of the theorem
and its proof, I will then include a short chapter on elementary
Ring Theory; it will be a collection of standard results that should
only serve as a
reference to the following chapter.\\

A slightly simplified version of the actual assertion will finally
be explained in Chapter (\ref{ch:theorem}), together with a detailed
proof. The next chapter will then explain how this simplified
version differs from the original claim made by Kayal \& Poonen, and
what changes might be made to the proof in the previous chapter in
order
to adapt it to the ``full version''.\\

Finally I will, in chapter (\ref{ch:running}), provide the reader
with some background about algorithms and computations, and give a
brief discussion about the running time of the algorithm.

\chapter{Elliptic Curves} \label{ch:ecs}
I will first of all state a few definitions and standard results
from the study of elliptic curves.  As some of the proofs require a
few technical lemmas that are not directly relevant to this
dissertation I will omit most of them; they are standard bookwork
and can be found e.g. in \cite{sil3} and \cite{flynn} (N.B. those
sources also provide the interested reader with a thorough insight
into Elliptic Curves).

\section{Preliminary Definitions}

Throughout these definitions, we shall denote by $K$ some field.

\newtheorem{predef}{Definition}[section]
\begin{predef}
$A_n(K) = \{(x_1,\dots,x_n): x_1,\dots,x_n \in K\}$, is called
\textit{affine n-space}.\end{predef} \begin{predef} When $P \in
A_n(K)$, we say that $P$ is \textit{K-rational} or \textit{defined
over K}. \end{predef} \begin{predef} Let $P_n(K):=
\{(x_0,\dots,x_n): x_0,\dots,x_n \in K, \textrm{not all 0}\}$,
subject to the relation that $(x_0,\dots,x_n) = (y_0,\dots,y_n) \in
P_n(K)$ if there exists $r \in K$, $r \ne 0$, such that
$(y_0,\dots,y_n) = (rx_0,\dots,rx_n)$. $P_n(K)$ is called
\textit{projective n-space over K}. \end{predef} \begin{predef} A
polynomial in $n$ projective variables is an \textit{(n +
1)-variable homogeneous polynomial}.\end{predef} \begin{predef} A
\textit{projective curve in $P_2$} is defined by a homogeneous
polynomial in 3 variables $F(X, Y, Z) = 0$. \end{predef}
\begin{predef} Let $C: f(x,y) = 0$ be an affine curve and let $P =
(x_0,y_0)$ be a point on C.  We say that P is a \textit{singular
point on C} if
\begin{displaymath} \frac{d}{dx} (f) |_P = \frac{d}{dy} (f) |_P = 0.\end{displaymath}
\item  A curve is called \textit{non-singular} if it does not contain any
singular points.\end{predef}

Finally, we are in a position to unambiguously define elliptic curves:\\

\begin{predef}[Elliptic Curves]  An \textit{elliptic curve over a field K} is a
non-singular, projective cubic curve, defined over $K$, with a
$K$-rational point. \end{predef}

\begin{predef} Let $C: f(x, y) = 0$ and $C': f(x, y) = 0$ be curves over $K$. A
\textit{rational map $\phi$} over $K$ from $C$ to $C'$ is a map
given by a pair $\phi_1, \phi_2$ of rational functions in $(x, y)$,
defined over $K$, with the property that given any point $P = (x_0,
y_0)$ on $C$, then $(\phi_1(x_0, y_0), \phi_2(x_0, y_0))$ lies on
$C'$. \\
If there also exists a rational map $\psi$ from $C'$ to $C$ such
that $\psi\cdot\phi$ is the identity on $C$ and $\phi\cdot\psi$ is
the identity on $C'$ then we say that $\phi$ is a \textit{birational
transformation} over $K$ from $C$ to $C'$, and that $C$ and $C'$ are
\textit{birationally equivalent} over $K$.

\end{predef}

\paragraph{Remark about the terminology}

An elliptic curve is not to be confused with an ellipse, which is a
plane algebraic curve usually given in the form \begin{displaymath}
\frac{y^2}{a^2} + \frac{x^2}{b^2} = 1 \end{displaymath} for some
non-zero constants $a, b$ in some field. There is however an
explanation for
the terminology.\\

Consider the relationship between the trigonometric functions sine,
cosine and tangent, and the arc lengths of a circle.  The further
study of elliptic curves shows that there is a similar relationship
between elliptic curves and arc lengths on ellipses.  These give
rise to so-called elliptic integrals of the form

\begin{equation} \label{eq:I} \int \frac{dx}{4x^3 + Ax + B}
\end{equation}\\

Integrals like (\ref{eq:I}) are multi-valued and only well-defined
modulo a period lattice $L$.  The "inverse" function of those
integrals
is a doubly periodic function called an elliptic function.\\

In fact every such function $P$ with periods independent over $L$
satisfies an equation of the form \begin{equation} P'^2 = 4 P^3 + AP
+ B \end{equation}

If we consider $(P, P')$ as a point in space then we can define a
mapping from the solutions of this equation to the curve

\begin{equation} Y^2 = X^3 + AX + B \end{equation}

This is the standard form for an elliptic curve that we shall be
concerned with throughout this dissertation.

\section{Arithmetic on Elliptic Curves}

\begin{predef}[Addition and Inverses] Let $C$ be an elliptic curve over a field $K$.  Let $o$ be its
K-rational point.  For any two points $a$, $b$ on $C$, denote by
$l_{a,b}$ the line through $a$ and $b$; if $a=b$ then $l_{a,b}$ is
defined to be the tangent to $C$ at $a = b$.  Let $d$ be the third
point of intersection of $l_{a,b}$ with $C$. Define $c$ to be the
third point of intersection between $C$ and $l_{o,d}$, the line
through $o$ and $d$. We then define $a + b := c$.\\

Let $k$ be the third point of intersection between $C$ and
$l_{o,o}$, the tangent to $C$ at $o$.  Let $a'$ be the third point
of intersection between $C$ and $l_{a,k}$.  Define $-a:= a'$.
\end{predef}

\subsubsection{Comment}

It is often convenient to write elliptic curves in affine form,
although it should be understood that we always mean a projective
curve.  For example, $C: y^2 = x^3 + 1$ will be used as the
shorthand notation for the projective curve $C: ZY^2 = X^3 + Z^3$.\\

It can be shown that any elliptic curve over K can be birationally
transformed over K to the Weierstrass form

\begin{equation}\label{eq:complicated} C: y^2 + a_1xy + a_3y = x^3 + a_2x^2 + a_4x + a_6
\end{equation}

To further simplify the equation we can use the following theorem.

\newtheorem{theorem}[predef]{Theorem}
\begin{theorem}
Let $K$ be a field with $char(K) \ne 2$.  Then any elliptic curve
over K is birationally equivalent over $K$ to a curve of the form
$Y^2 =$ cubic in $x$.  If $char(K) \ne 2 \textrm{ or } 3$ then we
can further reduce (\ref{eq:complicated}) to the form
\begin{equation} \label{eq:standard} Y^2 = X^3 + AX + B.
\end{equation}
\end{theorem}

For the purpose of this dissertation, we shall only be concerned
with elliptic curves over fields of characteristic $q > 3$, hence
(\ref{eq:standard}) will be treated as our standard equation for an
elliptic curve.  We will also adapt the convention to choose $o =
(0,
1, 0)$, the point at infinity.\\

Note that $Z = 0$ meets $C$ at $o$ three times.  Given $a = (x, y,
z)$, the third point of intersection between the curve and the line
through $a$ and $o$ is $(x, -y, z)$, which must then be $-a$. This
leads to the following result.

\newtheorem{Lemma}[predef]{Lemma}
\begin{Lemma} \label{def:gp} For an elliptic curve $C$ written in our
standard form (\ref{eq:standard}), we can simplify the formulae for
addition and inverses of points on $C$ as follows:
\begin{itemize}
\item  $-(x, y) := (x, -y)$ \\
\item If $d:= (x_3, y_3)$, the third point of intersection of $C$ and $l_{a,b}$,
then $a + b = (x_3, -y_3)$.
\end{itemize}
\end{Lemma}

\section{The Group Structure}

Now let us have a closer look at the rational points on an elliptic curve $C$.
With the above definitions of addition and $o$ we can show the following:\\

After a few computations it is easy to see that for all $a$ and $b$
on $C$ we have

\begin{itemize}
\item $a + b = b + a$ , \\
\item $a + o = o + a = a$ , and\\
\item $a + (-a) = (-a) + a = o$ .
\end{itemize}

Moreover, further computations that involve a few technical lemmas, will reveal that for $c$ on $C$ we also have
\begin{itemize}
\item $(a + b) + c = a + (b + c) \, .$
\end{itemize}

From this we can deduce the following theorem:

\begin{theorem}
Let $C$ be an elliptic curve over $K$.  The points on $C$, together
with the operation $a + b$ as defined in Lemma \ref{def:gp}, form a
group. The point $o$ acts as the identity in this group, and
inverses are given by $-a$ as in the definition above.
\end{theorem}

For a natural number $m \in \mathbb{N}$ we will from now on adapt the notation $[m]P:= P + P + \dots + P$ ($m$ times).
This map is also known as the ``multiplication-by-$m$-map'' from the curve to itself.
We can extend the definition of this to $m \in \mathbb{Z}$ by defining $[o]P:= 0$ and
$[-m]P:= -[m]P$.\\

So for example, if we have $P = (0, 1)$ on $C: Y^2 = X^3 +1$, then
$-[2]P = -(P+P)$.  Computing the tangent at $P$, we obtain the line
$L: Y = 1$, and the ``third point
of intersection'' of $C$ and $L$ being again $P$.  So $[2]P = -P = (0, -1)$, and hence $-[2]P =P$.\\

This map plays a central part in elliptic curve cryptography; its applications will later on be
extremely useful in this dissertation.

\section{Elliptic Curves over Finite Fields}

Now let us consider an elliptic curve $C$ over a finite field
$\mathbb{F}_q$.  Recall that we are only considering fields of
characteristic $q>3$ here.  The cases for $q = 2$ or $3$ are
similar, and some of our computations and notations could be adapted
to include those cases, too.  However, for the entire purpose of
this dissertation, those two cases will be irrelevant and we will
therefore exclude them in all our computations.\\

Consider the group of rational points on $C$ over $\mathbb{F}_q$.
The following result should be immediately obvious:

\begin{theorem}
Over a finite field $\mathbb{F}_q$, the number of rational points on
an elliptic curve $C$ is finite. \end{theorem}

We shall denote this number by $\#C(\mathbb{F}_q)$.  It may be asked
whether we can find out anything about this quantity.  The answer to
this is that we can indeed, and in fact the computation of this
number lies right at the heart of the
proof of the theorem.  \\

Before I give an in-depth discussion of how to compute the actual
value of $\#C(\mathbb{F}_q)$, I will give an upper and lower bound
on it, and define a few tools that we will later on need in our
computations.

\subsubsection{Discussion}

Consider a straight line $L$ over $\mathbb{F}_q$, given by $L: Y =
aX + b$.  What do we know about the number of points on $L$? \\

For every possible value $x \in \mathbb{F}_q$, i.e. $x \in
{0,1,\dots,q-1}$, there exists exactly one solution in
$\mathbb{F}_q$ for $y$, so we obtain $q$ rational points.  Also, the
point at infinity is always a rational point; in total we therefore
have exactly $q + 1$
rational points on $L$.\\

Now we can consider the  number of points on a curve $C$ of the form
$Y^2 = f(X)$ in a similar way: for each of the $q$ possible values
for $X$ we have one of the three cases:\\
\begin{itemize}
\item If $f(x)$ is a quadratic residue modulo $q$, we obtain
\textit{two} solutions for $Y$, namely $y = \pm \sqrt{f(x)}$;
\item If $f(x)$ is a quadratic non-residue modulo $q$, we will have
\textit{no} solutions for Y;
\item If $f(x) = 0$ we have precisely \textit{one} solution
for Y, namely $y = 0$.\\
\end{itemize}

From elementary Number Theory we know that in $\mathbb{F}_q$,
exactly half the values in $\{1,\dots,q-1\}$ are quadratic residues,
so we would expect the number of rational points on $C$ to be
roughly $q$ to represent the ``fifty-fifty chance of $f(x)$ being a
quadratic residue'', hence yielding $2$ solutions.  We then add the
point at infinity, and obtain as a rough estimate $q + 1$ rational
points.

\begin{predef}[Trace of Frobenius]

For a given curve $C$ over $\mathbb{F}_q$, the \textit{trace of
Frobenius $t$} is the quantity defined by the relation
$\#C(\mathbb{F}_q) = q + 1 - t \, .$\end{predef}

$t$ can therefore be regarded as the``error term'' in our estimate
of $\#C(\mathbb{F}_q)$.  The following theorem gives a bound on this
error term:

\begin{theorem}[Hasse's Theorem]\label{hasse}

\begin{equation} | t | \leq 2 \sqrt{q} \end{equation} \end{theorem}

A detailed proof of this can be found in \cite{sil1}. In Subsection
\ref{subs:zeta} of Chapter \ref{ch:orig} in this paper we will
see an alternative argument to deduce this.\\

A map that should be well-known to anyone who has studied basic
algebra and number theory is the Frobenius map.  It has a very
interesting property that has important applications in the study of
elliptic curves, as we shall soon see.

\begin{predef}

The $q^{th}$-power Frobenius map $\phi$ defined on an elliptic curve
$C$ over $\mathbb{F}_q$ maps points on $C$ to points on $C$ as
follows:

\begin{displaymath}
\phi = \left\{ \begin{array}{ll} C(\mathbb{F}_q) \to &
C(\mathbb{F}_q) \\
(x, y) \, \mapsto & (x^q, y^q)\\
 o \qquad \mapsto & o
\end{array} \right.
\end{displaymath}
\end{predef}

It is easily verified that $\phi$ is a group endomorphism for the
group of rational points on $C$ over $\mathbb{F}_q$ and is therefore
most commonly referred to as the Frobenius endomorphism.\\

As mentioned above, a deeper study of $\phi$ will reveal several
interesting results; the following property of $\phi$ is crucial to
this dissertation, and deserves particular attention.

\begin{Lemma}

The Frobenius endomorphism has characteristic polynomial\\

\begin{displaymath} \chi (X): X^2 - tX + q \, . \end{displaymath}
\end{Lemma}

\subsubsection{Outline Proof}
A full proof of this involves a lot of technical Lemmas; I will
therefore only state the main idea of the proof.\\
It relies on the fact that \begin{displaymath}\#Ker(\phi - 1) =
deg(\phi - 1) = q + 1 - t \, .\end{displaymath}

From this we can deduce that if we take an integer $m \ge 1$ with
$gcd(m, q) = 1$, then we have
\begin{displaymath}det(\phi_m) \equiv q \,(\textrm{mod }m), \qquad
tr(\phi_m) \equiv t\, (\textrm{mod }m) \, .\end{displaymath}

The details of this proof can be found in \cite{sil1}.\\

Clearly this is the same as writing

\begin{equation} \label{eq:frob} \phi ^2
- [t]\phi + [q] = [o].
\end{equation}

\newtheorem{cor}[predef]{Corollary}\label{cor:cru}
\begin{cor}
Hence we have that for a point $P:= (x, y)$ on $C$:
\begin{equation}
(x^{q^2}, y^{q^2}) - [t] (x^q, y^q) + [q] (x, y) = o.
\end{equation}

\end{cor}

\chapter{Counting Rational Points on Elliptic
Curves} \label{ch:counting}

As mentioned in Chapter (\ref{ch:intro}), the key to the proof of
our assertion is part of an algorithm that reveals the factors of
$h(z)$.\\

Although this dissertation is not about the efficient computation of
the number of rational points on elliptic curves, the algorithm that
we will later on adapt is in its original form a point-counting
algorithm for elliptic curves over finite fields. I will therefore
give a brief introduction to such
algorithms in general.\\

As mentioned earlier, over a finite field the number of rational
points on an elliptic curve is clearly finite.  In Chapter
(\ref{ch:ecs}) we
have seen an upper and lower bound for the number of such points.\\
Now it may be asked if we can actually compute the precise number of
rational points on a given curve. The answer is that we can indeed,
and there are several methods that can be applied to do this.\\

An explicit formula for the number of rational points on an elliptic
curve $C: Y^2 = X^3 + AX + B$ over $\mathbb{F}_q$ is given by the
following sum:
\begin{equation*} \#C(\mathbb{F}_q) = 1 + \sum_{x \textrm{ mod } q}
\Bigg( \bigg(\frac{X^3 + AX + B}{p} \bigg) + 1 \Bigg).
\end{equation*}

where $(\frac{a}{p})$ denotes the Legendre Symbol.\\

Computing $\#C(\mathbb{F}_q)$ this way takes $O(q^{1+\epsilon})$ bit
operations \footnote{For a definition of ``bit operations'', see
Chapter (\ref{ch:running}).}; this is clearly not very practical
when $q$ is a large prime number
(which, in practical applications of our theorem, it usually will be!).\\

Due to Rene Schoof however, we have a more efficient way of
computing $\#C(\mathbb{F}_q)$.  In his paper \cite{schoof}, Schoof
gives an explicit deterministic algorithm to compute the exact
number of points on any given curve over a finite field.\\

In the next chapter I will go into detail about this particular
algorithm, but first we will need yet more technical tools in order
to understand the
algorithm better.\\

As Corollary (\ref{cor:cru}) suggests, the algorithm involves
calculating coordinates of rational points of the form $[m]P$ where
$P$ is a rational point and $m$ a positive integer.  We will
therefore need techniques to efficiently compute those
coordinates.\\

It should be clear that the coordinates of $P_1 + P_2$ are rational
functions of the coordinates $(x_1, y_1)$ of $P_1$ and $(x_2, y_2)$
of $P_2$. By repetition of this calculation we can see that
multiplication by $[m]$ given by
\begin{displaymath}
(x,y) \mapsto [m](x, y)
\end{displaymath}

can also be expressed in terms of rational functions in $x$ and $y$.
Explicitly, we have the formulae given in the following section.

\section{The Division Polynomials}
\begin{Lemma} Let $C$ be an elliptic curve defined over
a field $K$ and let $m$ be a positive integer.  There exist
polynomials $\psi_m, \theta_m, \omega_m \in K[x, y]$ such that for
$P = (x, y) \in C(K)$ with $[m]P \ne o$ we have \begin{equation}
\label{eq:mP} [m]P = \bigg( \frac{\theta_m(x, y)}{\psi_m(x, y)^2},
\frac{\omega_m(x, y)}{\psi_m(x, y)^3} \bigg)
\end{equation}
\end{Lemma}

The polynomial $\psi_m$ is generally referred to as the
\textit{$m^{th}$ Division Polynomial} of $C$. $\theta_m$ and
$\omega_m$ can both be expressed in terms of $\psi_m$ as shown in
the explicit recursive expressions for $\psi_m$ below.

\subsubsection{Remark}
The expressions for $\psi_m$ that I will give are simplified for the
case where we can write our curve in the form $C: Y^2 = X^3 + AX +
B$.  They are given in a more general form in \cite{blake} for the
general curve

\begin{displaymath}
C: Y^2 + a_1 XY + a_3 Y = X^3 + a_2X^2 + a_4 X + a_6 \,
,\end{displaymath}
which could also be defined over fields of
characteristic 2 or 3.

\subsection{Explicit Expressions for $\psi$}\label{fdef}

Let $C: Y^2 = X^3 + AX + B$ be defined over $K$.\\

Then $\psi_m$ can be computed as follows:

\begin{eqnarray}
\lefteqn{\psi_0 = 0, \psi_1 = 1,  \psi_2 = 2y, {}} \nonumber \\
& & {} \psi_3 = 3x^4 + 6Ax^2 + 12Bx - A^2, {} \nonumber \\
& & {}\psi_4 =
4y(x^6 + 5Ax^4 + 20Bx^3 - 5A^2x^2 - 4ABx - 8B^2 - A^3), {} \\
& & {}\psi_{2m+1} =
\psi_{m+2} \psi_m^3 - \psi_{m-1} \psi_{m+1}^3, \, m \ge 2, {} \nonumber \\
& & {}\psi_{2m} = \frac{(\psi_{m+2} \psi_{m-1}^2 - \psi_{m-2}
\psi_{m+1}) \psi_m}{2y}, \, m > 2 {} \nonumber
 \end{eqnarray}

We can now in turn define $\theta_m$ and $\omega_m$ in terms of the
division polynomials:
\begin{eqnarray*}
\lefteqn{\theta_m = x \psi_m^2 - \psi_{m-1} \psi_{m+1} {}}\\
& & {}\omega_m = \frac{\psi_{2m}}{2 \psi_m}
 \end{eqnarray*}

Finally, we define \begin{equation} \label{eq:f} f_m := \left\{ \begin{array}{ll}
\psi_m & m\textrm{ odd} \\
\psi_m / (2y) & m\textrm{ even} \end{array} \right.
\end{equation}

The proof of these formulae involves straightforward but lengthy
calculations and will therefore be omitted; some more detail is
included in \cite{lang1}. It is however important to note the
following two facts:

\begin{cor} Let $f_m$ be defined as in (\ref{eq:f}).  Then

\begin{enumerate}
\item $f_m$ is a polynomial in $x$ only.
\item \label{degf} The degree of $f_m$ is at most $(m^2 -1)/2$ if $m$ is
odd, \\and at most $(m^2 - 4)/2$ if $m$ is even. \end{enumerate}
\end{cor}

The latter fact will be relevant in Chapter \ref{ch:running} when
calculating the running time of the algorithm.

\section{The $m$-Torsion Subgroup}

Clearly, when $K$ is a finite field, $\mathbb{F}_q$ for some prime
$q$, then $C(K)$ is a torsion group, i.e. every point on the curve
has finite order (since $C(K)$ itself is finite).  For a
non-negative integer $m$, the set of \textit{m-Torsion points} on
$C(K)$ is defined by \begin{displaymath} C[m]:= \{P \in C(K) | [m]P
= o\}.\end{displaymath}

We can now also express Corollary \ref{cor:cru} in terms of elements
of this subgroup:

\begin{cor}

For points of order $m$ on $C$, i.e. $P = (x,y) \in C[m]$, we have

\begin{equation}
\phi_m^2(P) - [\tau]\phi_m(P) + [k](P) = [o](P) = o \, ,
\end{equation}

where we define $\tau \equiv t$ (mod $m$) and $k \equiv q$ (mod
$m$).\end{cor}

It is easily verified that this is a subgroup of $C(K)$.  By
definition, $o \in C[m]$. The $m^{th}$ division polynomial $\psi_m$
characterises the $m$-Torsion subgroup as stated in the following
theorem.

\begin{theorem} Let $P \in C(K) \backslash {o}$, and let $m \ge 1$.
Then
\begin{equation*} P \in C[m] \Leftrightarrow \psi_m(P)
= 0. \end{equation*}

\end{theorem}

Clearly, this condition is equivalent to the following corollary,
which is more useful for our computations later:

\begin{cor} Let $P = (x, y) \in C(K) \backslash \{o\}$ be such
that $[2]P \ne o$ and let $m \ge 2$.  Then
\begin{equation}\label{le:tor}  P \in C[m] \Leftrightarrow f_m(x) =
0.
\end{equation}
\end{cor}

The 2-torsion points are excluded in this, since they satisfy
$\psi_2{P} = 0$, which we need to divide $\psi_m$ by in order to
obtain $f_m$ if $m$ is even.  However, we can immediately recognise
points of order 2 due to the fact that their $y$-coordinate is
always equal to zero (the reader may check this as an easy exercise
to become familiar with the arithmetic on elliptic curves).\\

To finish this rather technical section off, I will give an explicit
expression for $[m]P$, which is again just a straightforward
transformation of (\ref{eq:mP}):

\begin{equation} \label{eq:mP2}
[m]P = \bigg(x - \frac{\psi_{m-1} \psi_{m+1}}{\psi_m^2},
\frac{\psi_{m+2} \psi_{m-1}^2 - \psi_{m-2} \psi_{m+1}^2}{4y
\psi_m^3}\bigg)
\end{equation}

In the actual application of this result, $\psi_m$ will be replaced
by $f_m$ so that $[m]P$ is a rational function of $x$ only.  We will
show this explicitly later.\\

In the following chapter I will give a detailed explanation of the
deterministic algorithm that provides us with an efficient method to
count the number of rational points on a given elliptic curve over a
finite field.

\chapter{Schoof's Algorithm} \label{ch:schoof}

Schoof's Algorithm (published in April 1985) provides us with a tool
to compute $\#C(K)$, where $C$ is an elliptic curve given in
Weierstrass form, and $K = \mathbb{F}_q$ is a finite field.  The
algorithm takes $O((\log q)^8)$ elementary operations  and is
deterministic \footnote{Note that in his paper, Schoof shows that
the running time of his algorithm is $O((\log q)^9)$; it can however
be shown that one can make improvements on this bound.  This will be
discussed in more depth in Chapter \ref{ch:running}.}; it does not
depend on any unproved hypotheses. As usual, we will restrict
ourselves to the case where $char(K) \ne 2$ or $3$; those cases
again need separate treatment, however as mentioned before, they are
irrelevant for our
purposes.\\

I will first of all list the main steps of the algorithm, so that
the reader can refer to them when working through the following
section. Note that the purpose of some of those steps may not seem
immediately obvious, and some notation may be unfamiliar, but the
details will of course be filled in afterwards.

\section{Outline}

\texttt{INPUT: An elliptic curve $C: Y^2 = X^3 + AX + B$ defined
over $\mathbb{F}_q$ where $q$ is a prime $\ne 2$ or $3$.}
\texttt{
\begin{enumerate}
\item Compute the quantity $l_{max}$ defined in (\ref{eq:defl}) in Section
\ref{expl}
\item \label{step2} Set $l = 3$. Compute $\tau = t$ (mod
$l$) as follows:
\begin{enumerate}
\item \label{begin} Set $\tau \equiv 0$ (mod $l$).\\
\item \label{test} To test whether there exists a point $P$ on
$C$ such that \\ $\phi_l^2P + [k]P = \pm [\tau]\phi_l P$, where $k \equiv q$ (mod $l$):\\
Compute $H_{k,\tau}$ as defined in (\ref{eq:H0}) and (\ref{eq:Htau}) in Section \ref{expl}.\\
\item \label{gcd} Compute $gcd(H_{k, \tau}, f_l)$, where $f_l$ is as defined in (\ref{eq:f}) in Section \ref{fdef}.
\begin{enumerate}
\item If gcd$(H_{k, \tau}, f_l) = 1$, go to Step (\ref{nextt}).\\
\item If gcd$(H_{k, \tau}, f_l) \ne 1$, determine the correct sign of $\tau$
using the methods explained in Section \ref{expl}.  Set $\tau = \pm \tau$ accordingly.  Go to step (\ref{output}).\\
\end{enumerate}
\item \label{nextt} Set $\tau = \tau + 1$.  Go to step (\ref{test}).\\
\item \label{output} OUTPUT: ($\tau, l$).  Set $l =$ next prime $\le l_{max}$, go to step (\ref{begin}).
If the next prime $l > l_{max}$, go to step (\ref{crt}).\\
\end{enumerate}
\item \label{crt} Compute $t$ using the Chinese Remainder Theorem applied to
$(\tau, l)$ for all $l$.\\
\item  Compute $\#C(\mathbb{F}_q)$.\\
\item FINAL OUTPUT: $\#C(\mathbb{F}_q)$.
\end{enumerate}
}
\section{Explanation}\label{expl}

Now the above looks very abstract and clearly requires
explanation.\\

Note that only in the last step we are concerned with
$\#C(\mathbb{F}_q)$. The algorithm actually computes the trace of
Frobenius; this is clearly equivalent to computing
$\#C(\mathbb{F}_q)$ due to the one-to-one correspondence between the
two quantities, $\#C(\mathbb{F}_q) = q + 1 - t$.\\

Since Hasse's Theorem \ref{hasse} provides us with a bound on $t$,
it will be sufficient to compute $t$ modulo a sufficiently large
number of primes and then recover the value of $t$ by an application
of the Chinese Remainder Theorem.\\

I will explain the algorithm for $\tau \equiv 0$ in thorough detail
first; the case for $\tau \in \{1,\dots,(l-1)/2\}$ will then only be
outlined.  It follows the same idea, but involves computing
slightly more complicated polynomials.\\

Since the algorithm is in its abstract form very technical and hence
somewhat difficult to follow, I will, as an example of its
application, demonstrate each step by performing it on the elliptic
curve $C: Y^2 = X^3 +1$ over the
finite field $\mathbb{F}_5$.\\

Note that this example is almost trivial, since over $\mathbb{F}_5$
we can find the number of points by inspection rather easily.

The example will however also show that the computations over such
small fields already involve very complicated looking polynomials.
In practice, we would apply Schoof's Algorithm to fields of
characteristic a large prime. \footnote{In order to avoid confusion
I will print the example in blue so that the reader can easily
distinguish more easily between ``theory'' and ``practice'', since I
will often skip between the two.}\\

We begin by defining $l_{max}$ to be the smallest prime such that
\begin{equation} \label{eq:defl} \prod_{\substack{l \textrm{ prime} \\ 2\le l \le l_{max}}} l
> 4 \cdot \sqrt{q} \, . \end{equation}

This bound is sufficient for us to obtain enough values of $\tau$
(mod $l$) to recover $t$ using the Chinese Remainder
Theorem.\\

We know from (\ref{eq:frob}) that the trace of Frobenius $t$
satisfies \begin{displaymath} \phi ^2 + q = t \phi,
\end{displaymath} so if we reduce this equation modulo $l$ we have
\begin{equation} \label{eq:tau} \phi _l ^2 + k = \tau \phi_l \textrm{
, where } \tau \equiv t \textrm{ mod }l, \, k \equiv q \textrm{ mod
} l\end{equation}

for all points on $C$ of order $l$, i.e. for $P = (x, y) \in C[l]$.
\newpage

In order to find $\tau$ we need to check for which $\tau \in
\{0,1,\dots,l\}$ the relation (\ref{eq:tau}) holds.  To do this, we will test for
$\tau \in \{0,1,\dots,(l-1)/2\}$
whether a point $P = (x, y)$ exists in $C[l]$ such that
\begin{equation} \label{eq:t} \phi_l^2 P + [k]P = \pm [\tau] \phi P \, .\end{equation}

By applying (\ref{eq:mP2}), we can see that this is equivalent to
testing for which $\tau$ we have
\begin{eqnarray} \label{eq:expl} (x^{q^2}, y^{q^2}) + \Big(x - \frac{\psi_{q-1}
\psi_{q+1}}{\psi_q^2}, \frac{\psi_{q+2} \psi_{q-1}^2 - \psi_{q-2}
\psi_{q+1}^2}{4y \psi_q^3}\Big) \nonumber \\
\nonumber \\
=  \left\{
\begin{array}{ll} 0 &
\textrm{if } \tau \equiv 0 \textrm{ (mod l)}\\
\bigg(x^q - \Big(\frac{\psi_{\tau - 1} \psi_{\tau + 1}}{\psi_\tau
^2}\Big)^q, \Big(\frac{\psi_{\tau + 2} \psi_{\tau - 1} ^2 -
\psi_{\tau - 2} \psi_{\tau + 1}^2}{4y \psi_\tau ^3} \Big)^q \bigg) &
otherwise. \end{array} \right. \end{eqnarray}

So let us now run through the algorithm to see what happens at each
step. First of all we compute $l_{max}$.

\subsection{Step (\ref{step2})}

We set $l = 3, \tau = 0$. We then test whether there exists a point
in $C[l]$ such that
$\phi_l^2P = \pm [k]P$.\\

Comparing the x-coordinates of both sides in (\ref{eq:expl}), we can
see
 that this holds if and only if
\begin{equation} x^{q^2} = x - \frac{\psi_{k-1}\psi_{k+1}}{\psi_k^2}
\end{equation}.

In order to obtain a univariate polynomial in $x$ only, we replace the $\psi_n$ by $f_n$
and multiply through by the denominator.  Let us now define $H_{k,0}(x)$ as follows:
\begin{equation} \label{eq:H0} H_{k,0}(x) : = \left\{
\begin{array}{ll}
(x^{q^2} - x) f_k^2(x)(x^3 + Ax + B) + f_{k-1}(x)f_{k+1}(x) &
 (k \textrm{ even)} \\
(x^{q^2} - x)f_k^2(x) + f_{k-1}(x)f_{k+1}(x)(x^3 + Ax + B) &
 (k \textrm{ odd).}
\end{array} \right. \end{equation}

\subsection{Step (\ref{gcd})}

We have now reduced the problem of testing whether relation
(\ref{eq:tau}) holds for $\tau \equiv 0$ (mod $l$) to checking
whether there exists a $P = (x, y) \in C[l]$
such that $H_{k,0}(x) = 0$.\\

\textcolor{blue}{Let us consider our example: Since we have $q = 5$,
we can take $l_{max} = 5$. We will now compute $f_n$ now for n =
0,\dots,4:} \textcolor{blue}{\begin{eqnarray*}\label{fsex}
\lefteqn{f_0(x) = 0, f_1(x) = 1, f_2(x) = 1, {}}\\
\\
& & {} f_3(x) = 3x^4 + 12x, f_4 = x^6 +20x^3 -8.{}
\end{eqnarray*}}
\textcolor{blue}{Since we have set $l = 3, \tau = 0$, now need to
compute $H_{2,0}$.}
\textcolor{blue}{
\begin{eqnarray}
\lefteqn{H_{2,0} = (x^{25}-x)(x^3+1)f_2(x)^2 + f_1(x)f_3(x) {}} \\
\nonumber \\
& & {} = (x^{25} - x)(x^3+1) + f_3(x) {} \nonumber
\end{eqnarray}}

Although it may not seem immediately obvious, why this is any easier
than the original problem, it is indeed a simplification: We recall
from (\ref{le:tor}) that \begin{displaymath} P = (x, y) \in C[l]
\Leftrightarrow [l]P = 0 \Leftrightarrow f_l(x) = 0 \,
.\end{displaymath} On the other hand we know that if our chosen
$\tau$ is indeed the trace of Frobenius,
then for all such $x$, we have $H_{k, \tau}(x) = 0$.\\

From this we deduce that all roots of $f_l$ are also roots of $H_{k, \tau}$, and
the two polynomials therefore have a non-trivial greatest common divisor.\\

So rather than attempting to solve the equation $H_{k,0}(x) = 0$, we
only need to compute the greatest common divisor of $H_{k,\tau}$ and
$f_l$; this explains step (\ref{gcd}) of the algorithm.  Note here
that in order to compute the greatest common divisor we use the
Euclidean Algorithm, which the reader should
be familiar with; it is briefly outlined in the next chapter.\\

Now consider the case where $gcd(H_{k,0}, f_l) = 1$.  This happens
if and only if $H_{k,0}$ and $f_l$ have no roots in common.  In that
case we have $H_{k,0}(x) \ne 0 \textrm{ for any }P = (x,y) \in C[l]$
and so we conclude that there exists no $P \in C[l]$ such that
(\ref{eq:tau}) holds. Clearly this means that $t \ne \tau$ (mod
$l$), and we go to step (\ref{nextt}) to set $\tau = \tau + 1$ and
try again for this new value of $\tau$.\footnote{Note here that we
never hit $\tau + 1 > (l-1)/2$ as we know that \textit{exactly} one
$\tau \in \{0,\dots,(l-1)/2\}$ will satisfy (\ref{eq:tau}) and so
$(l-1)/2$ is the largest value that $|\tau|$ can take. Once we have
hit this value we know it is the correct solution and we find
ourselves in Step (\ref{gcd}).ii, from which we proceed to
(\ref{nextt}) straight away.}  I will
return to this case later.\\

If, on the other hand, the greatest common divisor is non-trivial,
then we know that if we have a point $P$ in $C[l]$, it will
necessarily satisfy the desired
property (\ref{eq:tau}).\\

\textcolor{blue}{In our example, this step boils down to finding the
greatest common divisor of $(x^{25} - x)(x^3+1) + f_3(x)$ and
$f_3(x)$.  This in turn is just \begin{displaymath}gcd((x^{25} -
x)(x^3+1), f_3(x)) = gcd((x^{25} - x)(x^3+1), 3x^4 +
12x).\end{displaymath}} \textcolor{blue}{We find that this greatest
common divisor turns out to be $x$ and hence proceed to sub-step ii,
which is explained
below:}\\

\subsubsection{gcd $\ne 1, \, \tau =0$}

Now there are two ``subcases'' to be considered: Namely
when $\phi_l^2 P = [-q]P$ and when $\phi_l^2 P = [+q]P$.  We now run through a ``sub-algorithm''
for the case $\tau = 0$.

\begin{itemize}
\item Test whether $\phi_l^2P = - [k]P$
or $\phi_l^2P = + [k]P$ by checking the $y$-coordinate in a similar way.\\
\begin{itemize}
\item If $\phi_l^2P = - [k]P$, go to step (\ref{output}) of the main algorithm.\\
\item If $\phi_l^2P = + [k]P$, test whether $q$ is a square modulo $l$\\
\begin{itemize}
\item If $\big(\frac{q}{l}\big) = - 1$, go to step (\ref{output}) of the main algorithm.\\
\item If $\big(\frac{q}{l}\big) = 1$, let $\omega^2 = q$ (mod $l$) and
test whether $\phi_l P = [\omega]P$ or $\phi_l P = -[\omega]P$, .
Set $\omega_0 = \pm \omega$ accordingly.  Set $\tau = 2 \omega_0$,
go to step (\ref{output})
of the main algorithm.\\
\end{itemize}
\end{itemize}
\end{itemize}

Now why are we doing all this?\\

\paragraph{Case 1} First assume that $\phi_l^2 P = [-k]P$.  So we know that
$\tau \phi_l P =0$, and since $\phi_l P \ne 0$, we can conclude that
$t \equiv 0$ (mod $l$). So we proceed to Step (\ref{output}) in
the main algorithm and then run the algorithm for the next prime $l$.\\

\paragraph{Case 2} On the other hand, if $\phi_l^2 P = [+k]P$, then
\begin{displaymath} (2q - \tau \phi_l)P = 0 \qquad \textrm{and so} \qquad \phi_lP
= \frac{2q}{t}P. \end{displaymath} Let us apply $\phi_l$ to both
sides and use the equality satisfied by $P$; so we get
\begin{displaymath}qP = \phi_l^2P = \phi_l\big(\frac{2q}{t}P\big) = \big(\frac{2q}{t}\big)^2P \,
,\end{displaymath}

and hence that $t^2 \equiv 4q$ (mod $l$).  Again, we must split this
into two subcases:  When $q$ is a quadratic residue modulo $l$ and
when it is not.

\begin{itemize}
\item $\big(\frac{q}{l}\big) = - 1$: In this case we can
conclude that $\tau \equiv 0$ (mod $l$) and go to Step
(\ref{output}).
\item $\big(\frac{q}{l}\big) = 1$:  Let $0 < \omega < q-1$ denote a
square root of $q$ modulo $l$.  Since we have $(2q - \tau \phi_l)P =
0$, we can see that $2q/t$ is an eigenvalue of $\phi_l$; but $t/2 =
\pm \sqrt{q}$, so either $\sqrt{q}$ or $-\sqrt{q}$ is an eigenvalue
of $\phi_l$.  To test this, we proceed exactly as before in checking
whether $\phi_lP = [\pm \sqrt{q}]P$.  If we denote by $\omega$ the
correct eigenvalue, we can finally set $\tau = 2 \omega$ and proceed
to Step (\ref{output}).
\end{itemize}

\textcolor{blue}{In our example, we check the $y$-coordinates of
(\ref{eq:expl}) and proceed as above: assuming that $\phi_l^2P = -
[k]P$, we turn this into a polynomial in $x$ that depends on $k$ and
$\tau$, and compute its greatest common divisor with $f_3$.  We find
that this greatest common divisor is non-trivial and hence
$\phi_3^2P = - [2]P$ is indeed the correct solution.  So $t \equiv
0$ (mod $3$) in our example.}\\

\textcolor{blue}{Applying the same methods, we compute $\tau$ (mod
$5$)
and find that $t \equiv 0$ (mod $5$), too.}\\

\textcolor{blue}{Let us return to the algorithm to see what would
have happened if the greatest common divisor had been trivial.}

\subsubsection{gcd $= 1$}

Here we have that for no point $P \in C[l]$, relation (\ref{eq:t})
is satisfied.  From this we conclude that $t \ne \tau$ (mod $l$) and
so we need to check whether the next value of $\tau$ is the trace of
Frobenius mod $l$, i.e. whether $\phi_l^2 P + [k] P = \pm
[\tau]\phi_l P$ for $P\in C[l]$. So we go back to step (\ref{test})
and compute $H_{k, \tau}$ for this new value of $\tau$. Referring to
(\ref{eq:expl}) again, we know that
\begin{equation} \label{eq:tx} \big(\phi_l^2P + kP\big)_X = x^{q^2} + x +
\frac{f_{k-1}f_{k+1}}{f_k^2} + \lambda^2 + \lambda ,
\end{equation}

where \begin{displaymath} \lambda =
\frac{(y^{q^2}+y+x)xf_{k}^3+f_{k-2}f_{k+1}^2
+(x^2+x+y)(f_{k-1}f_kf_{k+1})}{xf_k^3 (x+x^{q^2})+
xf_{k-1}f_kf_{k+1}}.
\end{displaymath}

On the other side we have that

\begin{equation} \big(\pm\tau\phi_lP\big)_X = x^q +
\bigg(\frac{f_{\tau+1}f_{\tau-1}}{f_\tau^2} \bigg)^q.
\end{equation}

Now we can, in a similar way as above, transform the equation by
reducing modulo the curve equation so that we have polynomials of
degree at most one in $y$, since we can substitute $(x^3 + Ax + B)^m$ for any $y^{2m}$.\\

We then obtain an equation of the form $a(x) - yb(x) =0$, hence $y =
a(x)/b(x)$ for some $a(x), b(x) \in \mathbb{F}_q(x)$.  Again
substituting for $y$ in the curve
equation we therefore finally get
\begin{displaymath}y^2 = (x^3 + Ax + B) = \Big(\frac{a(x)}{b(x)}\Big)^2 \,
.\end{displaymath}

So we define \begin{equation} \label{eq:Htau}
H_{k,\tau}:=a(x)^2 - (x^3 + Ax + B)b(x)^2,
\end{equation}

a polynomial in $x$ only.\\

Now we can proceed precisely as before:  We want to check whether
for $x$ such that $f_l(x) = 0$ we also have $H_{k,\tau} = 0$, i.e.
whether the roots of
$f_l$ are also roots of $H_{k,\tau}$.\\
So in Step (\ref{gcd}) we compute the greatest common divisor of
$H_{k,\tau}$ and $f_l$.\\

If the points on $C[l]$ do not satisfy the Frobenius relation and
hence the greatest common divisor is 1, we
conclude that this value of $\tau$ is also not the correct value.\\
We are therefore sent to Step (\ref{nextt}) to proceed to the next
possible value of $\tau$ and then return to to Step (\ref{test}),
where we run the same test for that new value. \\

Otherwise we have that $t \equiv \pm \tau (l)$ for our chosen
$\tau$. In this case we need to check which sign is correct: We
refer to (\ref{eq:expl}) again, this time comparing the
$y$-coordinates of both sides, and check in a similar
manner which is the correct sign.\\

\subsection{The final steps}

This way we eventually obtain enough values for $t$ (mod $l$) so
that we can finally proceed to Step (\ref{crt}) and apply the
Chinese Remainder Theorem to the pairs $(\tau, l)$.  Finally we can
calculate $\#C(\mathbb{F}_q)$, which completes the algorithm.\\

\textcolor{blue}{Applying the Chinese Remainder Theorem to our
example $C: Y^2 = X^3 + 1$ over $\mathbb{F}_5$, where we had that
$t\equiv 0$ (mod $3$) and $t\equiv 0$ (mod $5$), we can deduce that
$t = 0$. So the number of rational points on $C$ over
$\mathbb{F}_5$ is $5 + 1 - 0 = 6$.}\\

\textcolor{blue}{Since we have chosen such a simple example, we can
verify this result by inspection, i.e. by trying each value of $x \in \mathbb{F}_5$:}\\
\textcolor{blue}{On $C: Y^2 = X^3 +1$ over $\mathbb{F}_q$ we have
the following rational points: \begin{displaymath} o, (0, \pm1), (2,
\pm2), (4, 0).\end{displaymath} So we obtain the same result; there
are 6 rational points on $C$ over $\mathbb{F}_5$.}\\

The topic of point-counting algorithms, improvements of Schoof's
Algorithm and its applications is an extremely interesting and
wide-ranging one;  I refer the interested reader to \cite{blake} and
\cite{kobl1} for the further study of this subject.  Any deeper
discussion about this subject is however irrelevant to this
dissertation.

\chapter{Some Ring Theory} \label{ch:ring}

For the purpose of a clearer understanding of the proof that will
follow in the next chapter, we will need to recall some elementary
theory about rings and fields. \\

The following results should be known to the reader.  I will
therefore omit proofs to the assertions made; they should be
regarded as a list of results that the reader may refer to in some
steps of the proof of the theorem.\\

Details of proofs can be found in e.g. \cite{cohn1}, \cite{cohn2} or
\cite{hers}.

\section{Elliptic Curves defined over a Ring}

Let $h(z)$ in $\mathbb{F}_q[z]$ be a nonzero polynomial of degree 2
with distinct roots in $\mathbb{F}_q$, and consider a curve $C$ over
the ring $R: = \mathbb{F}_q[z]/(h(z))$.  We may
view $C$ as a pair of curves over $\mathbb{F}_q$ as follows: \\

$C$ has coefficients of the form $(r z + s)$, where $r, s \in
\mathbb{F}_q$.\\

Since $R \cong \mathbb{F}_q \times \mathbb{F}_q$, we can apply the
isomorphic map $(r z + s) \mapsto (r a + s, - r a + s)$, where $ \pm
a$ are the roots of $h(z)$ in $\mathbb{F}_q$, to the curve to obtain
a pair of curves, both defined over $\mathbb{F}_q$. I.e.:
\begin{equation} \label{eq:C} C: Y^2 = X^3 + (\alpha z + \beta)X +
(\gamma z + \delta) \mapsto \left\{ \begin{array}{ll} C_{+}: & Y^2 =
X^3 + (\alpha
a + \beta)X + (\gamma a + \delta) \\
C_{-}: & Y^2 = X^3 + (-\alpha a + \beta)X + (-\gamma a + \delta)
\end{array} \right.
\end{equation}\\

\section{Rings and Fields}

The fundamental difference between a ring and a field is that in a
ring we may have non-units.  That is, we may have elements $r \in
R$, such that there exists no $s \in R$ with $r \cdot s = 1_R$.\\

A \textit{zero divisor} is an element $r \in R, r\ne 0$ such that
there exists $s \in R \backslash\{0\}$ with $r\cdot s = 0_R$.

\begin{Lemma} If $r \in R$ is a zero divisor, then it is a non-unit.
\end{Lemma}

\paragraph{Example}
For instance, in the ring $A = \mathbb{Z} / 4\mathbb{Z}$, we have
that $2 \ne 0$, but $2 \cdot 2 = 4 = 0_A$. 2 is therefore a zero
divisor.  It is also a non-unit: there is no element $s \in A$ such
that $s \cdot 2 = 1_A$.\\

Although this should be obvious, the following result is worth some
particular attention:

\begin{cor}\label{nonunits} In $R \backslash \{0\} =
\big(\mathbb{F}_q[z] / (z^2 - a^2)\big) \, \backslash \{0\}$, the
non-units are $(z - a)$ and $(z + a)$.
\end{cor}

\section{Euclid's Polynomial Division Algorithm}

Let us consider two univariate polynomials $f(x), g(x)$ defined over
some field $F$.  Euclid's Polynomial Division Algorithm provides us
with an efficient tool to compute the greatest common divisor of $f$
and $g$.

\subsection{Long Division of Polynomials}

Recall from school how we divide polynomials:  First we divide the
leading term of the higher degree polynomial by the leading term of
the lower degree polynomial.  Now think about what ``dividing''
means:  we try to find an element $a$ such that $a \cdot b = c$,
where $b$ is the leading coefficient of the lower degree polynomial
and $c$ that of the higher degree polynomial. All this should of
course be clear, but as it will be crucial later on, it is again
worth noting down the following result:\\

\begin{Lemma}\label{division}
We have $a = c\cdot b^{-1}$.\\
\end{Lemma}

Hence in order to find $a$, we compute the inverse of $b$ and
premultiply it by $c$.

\subsection{Euclid's Algorithm}

This is just a brief outline of the algorithm.  Details can be found
in any undergraduate book on linear algebra, e.g. \cite{cohn1},
\cite{cohn2} or \cite{hers}.

\subsubsection{Proposition}

For $f(x), g(x) \in F[X]$ with $0<deg(g)<deg(f)$, there exist
$r_1(x), q_0(x) \in F[X]$
such that we can write\\

$f(x) = g(x)q_0(x) + r_1(x)$, with $deg(r_1) < deg(f) \textrm{ or } r_1 \equiv 0$.\\

The polynomials $r_1(x)$ and $q_0(x)$ are computed by long division
of polynomials. As the next step in the Algorithm, we define a
sequence $r_i(x)$ as follows:\\

$r_0(x) = g(x)$\\

$r_i(x) = r_{i+1}(x)q_{i+1}(x) + r_{i+2}(x)  \textrm{ with } deg(r_{i+2})<deg(r_{i+1}) \textrm{ or } r_{i+2} = 0$ .\\

We eventually obtain\\

$r_{n-1}(x) = r_{n}(x)q_{n}(x) + r_{n+1}(x), \qquad r_{n+1}(x) = 0, r_{n}(x) \ne 0$ .\\

At this point the algorithm ends, and returns $r_{n}(x)$ as the greatest common divisor.\\

Now we are finally ready to tackle the actual problem we are aiming
to solve.

\chapter{The Theorem} \label{ch:theorem}

In this chapter I will discuss the theorem to be proved. First of
all I will give the already simplified version of the theorem and
prove it. The original statement of it is somewhat more complicated
and requires a few more definitions; this will be discussed in
Chapter \ref{ch:orig}.

\section{Statement of the Theorem}
\subsection{The problem}

The problem to be solved here is: \\

Find a deterministic algorithm with
\begin{itemize}
\item \texttt{INPUT: \begin{itemize}
\item A finite field $\mathbb{F}_q$ and
\item a nonzero polynomial $h(z)$ of degree 2 in $F_q[z].$
\end{itemize}
\begin{eqnarray}\label{eq:alg} \end{eqnarray}
\item OUTPUT: \begin{itemize}
\item The factors of $h(z)$ over $\mathbb{F}_q.$
\end{itemize}
\item Running time: polynomial in the size of the input, i.e., bounded by a
universal polynomial in $(1+ \textrm{deg }h)(\log q)$}.
\end{itemize}

\subsection{The hypothesis}

We are given a polynomial $h(z)$ of degree 2 with roots in the finite field
$\mathbb{F}_q$, where $q$ is a prime.\\

Assume that there exists an elliptic curve $C: Y^2 = X^3 + AX + B$ over the ring
$R: = \mathbb{F}_q[z]/(h(z))$ for which there exists a
prime $l$, such that we have $\#C_+(\mathbb{F}_q) \ne \#C_-(\mathbb{F}_q)$ (mod $l$).
\begin{equation} \label{qu:hyp}
\end{equation}

\begin{theorem}[Kayal]
Given (\ref{qu:hyp}), there exists an algorithm as in
(\ref{eq:alg}).
\end{theorem}

\section{The Proof}

Let us define $h(z):= (z^2 - a^2) \in \mathbb{F}_q[z]$ (where we do
not know $a$ but only $a^2$) and let
\begin{displaymath}C: Y^2 = X^3 + (\alpha z + \beta)X + (\gamma z + \delta)
\end{displaymath}

over the ring $R$ as above.  Let $l$ be the prime number that
satisfies the hypothesis of the theorem, i.e. such that
$\#C_{+a}(\mathbb{F}_q) \ne \#C_{-a}(\mathbb{F}_q)$ (mod $l$).

Let $t_+, t_-$ be the respective traces of Frobenius of $C_+$ and
$C_-$. Then we have that $t_+ \ne t_-$
(mod $l$).\\

The idea of the proof is that Schoof's point counting algorithm
\textit{is} an algorithm that solves our problem of factorising
$h(z)$.  I claimed earlier that when we apply it to
$C$, it will at some point reveal the factors of $h$.\\

Schoof's Algorithm is defined for elliptic curves over finite
fields, whereas $C$ is defined over a ring.  Note that if the
underlying hypothesis for the theorem were not fulfilled, we could
in general run the algorithm over curves defined over a
ring without any problems.\\

Running Schoof's Algorithm over $C$ is equivalent to running it over
$C_+$ and $C_-$ simultaneously.  Every operation that we are performing on $C$ can
be thought of as performing the same operations on $C_+$ and $C_-$ if we map $z$ to $\pm a$ accordingly.\\

Assume that we are in Step (\ref{begin}) of the algorithm with $l$,
the prime number with the desired property.  Checking every value of
$\tau \in \{0,\dots,(l-1)/2\}$ to see whether it satisfies $\phi^2 -
\tau \phi + q = 0$ is hence the same as checking whether there
exists a point $P^+ = (x^+, y^+)$ in $C_+[l]$ such that the relation
$\phi^2 - \tau \phi + q = 0$ is satisfied, and whether for a point
$P^- = (x^-, y^-)$ in $C_-[l]$, this equation holds.\\

Let $f_l^+$ and $f_l^-$ denote the $l^{th}$ division polynomial on $C_+$ and
$C_-$ respectively, and let $H^+_{k,\tau}$ and $H^-_{k,\tau}$ be as defined in (\ref{eq:Htau})
for the two curves accordingly.\\

Now, without loss of generality, we assume that $t_+ < t_-$.
Consider Step (\ref{gcd}) in the algorithm with $\tau = \tau_+
\equiv t_+$ (mod $l$). We compute $H^+_{k,\tau}$ and $H^-_{k,\tau}$.
Since $\tau = \tau_+$, we will
find that for all points $P^+ = (x^+, y^+)$ in $C_+[l]$, we have $H^+_{k,\tau}(x^+) = 0$.\\

On the other hand however, since $\tau \ne \tau_- \equiv t_-$ (mod $l$), we know
that \textit{no} point in $C_-[l]$ satisfies $H^-_{k,\tau} = 0$.\\

Now consider this step of the algorithm over $C$ itself. So we
attempt to compute $gcd(H_{k,\tau}, f_l)$, as usual, using the
Euclidean Algorithm.\\

Suppose we are trying to divide some polynomial $r_i(x)$ by
$r_{i+1}(x)$ where $deg(r_i) > deg(r_{i+1})$.  I will now make the
following claim:

\newtheorem{Proposition}[predef]{Proposition}
\begin{Proposition} \label{finalprop}

The leading coefficient $c(z)$ of $r_{i+1}(x)$ is a non-unit in $R$,
for some $i$.
\end{Proposition}

\begin{cor}\label{finalcor}
Proposition \ref{finalprop} completes the proof.
\end{cor}

\paragraph{Proof of Corollary \ref{finalcor}}

If Proposition \ref{finalprop} is true, then from Lemma
\ref{division}, we know that we are trying to compute the inverse of
$c(z)$.  Since this is a non-unit in $R$, it has no inverse, and
hence the
algorithm ``crashes''.\\

Now, the non-units in $R \backslash \{0\}$ are $(z \pm a)$ as noted
in Corollary \ref{nonunits} (possibly multiplied by a constant in
$\mathbb{F}_q)$.\\
So if we compute $gcd(c(z), h(z))$, we obtain a non-trivial factor
of $h(z)$, as required.\\

\paragraph{Proof of Proposition \ref{finalprop}}

Imagine that the lower degree polynomial $r_{i+1}(x)$ \textit{never}
has leading coefficient a non-unit in $R$. The Euclidean Algorithm
will just run smoothly over the ring as if it were a field.\\

Now recall that everything we are doing with the curve $C$
is equivalent to performing the same operations on the pair of curves $C_+$ and $C_-$ simultaneously.\\

Let us once more consider in detail the relationship between the
polynomials $r_i$ for $C$ and the $r_i^+$ and $r_i^-$ for $C_+$ and
$C_-$.  The latter two are just
evaluations of the coefficients of $r_i$ at $z = a$ and $z = -a$ respectively.\\

So if we assume that the leading coefficient of $r_i$ is a unit for
every $i = 0,\dots,n$ (where we have that $r_{n+1} = 0)$,
then the leading coefficient of $r_i$ \textit{never vanishes} on $C_+$ and $C_-$.\\

This implies that the degree of the polynomials $r_i$ is the same as
the degree of the $r_i^+$ and $r_i^-$.  In particular, we note that
for all $i = 0,\dots,n$, the degrees of $r_i^+$ and $r_i^-$ are the same.\\

So finally, we conclude that the degree of $r_n^+$ is the same as
the degree of $r_n^-$.  But recall that $r_n^+$ and $r_n^-$ are the
greatest common
divisors of $H_{k,\tau}^+$ and $f_l^+$, and $H_{k,\tau}^-$ and $f_l^-$ respectively.\\

By assumption however, the greatest common divisor of $H_{k,\tau}^-$
and $f_l^-$ is 1, since $\tau \ne \tau_-$ (mod $l$), whereas that of
$H_{k,\tau}^+$ and $f_l^+$ is strictly non-trivial!\\

This is clearly a contradiction.\\

We can now see that at some point we must
encounter a non-unit as the leading coefficient of some $r_i$.\\

By Corollary \ref{finalcor}, this completes the proof.\\

\chapter{The original statement of the Theorem} \label{ch:orig}

\section{Polynomials of higher degree}

As mentioned earlier, the full assertion made by Dr Kayal is
slightly more advanced.  Instead of restricting himself to
polynomials of degree 2, he claimed that the assertion would hold
for any polynomial with distinct roots in a finite field.\\

On closer inspection, one can see that this is plausible, and that
in fact the proof will be very similar to the one given above.  One
needs to think of an elliptic curve over the ring
\begin{displaymath} R:= \mathbb{F}_q[z] / (h(z)) = \mathbb{F}_q[z] /
\big((z-\alpha_1)(z-\alpha_2)(\dots)(z-\alpha_n)\big)
\end{displaymath}

as a family of curves over $\mathbb{F}_q$ in the same way as above,
i.e. with $z$ evaluated at $\alpha_i$ for each ``subcurve'' $C_i$.\\

If we are then given a curve $C$ over this ring, such that for some
prime $l$ the number of rational points on $C_i$ is not congruent to
the number of rational points on $C_j$ modulo $l$, for some $i < j$,
the same problem as discussed above will
arise in Schoof's Algorithm.\\

Let $\tau_i$ be the trace of Frobenius of $C_i$ modulo $l$, which is
hence not equivalent to the trace of Frobenius of $C_j$.  Adapting a
similar notation as before, and using the same arguments, we can
deduce that $gcd(H_{k,\tau_i}^{(i)}, f_l^{(i)})$ is strictly
non-trivial,
whereas $gcd(H_{k,\tau_i}^{(j)}, f_l^{(j)}) = 1$ . \\

When computing $gcd(H_{k,\tau_i}, f_l)$, on $C$, the Euclidean
Algorithm will again break down in an attempt to compute the inverse
of a non-unit in $R$, which we will inevitably encounter as the
leading coefficient of some $r_i$.  The reason for this is precisely
the same as in the case for $deg(h) = 2$:  if this never happened,
then we would be able to conclude from this fact that the greatest
common divisors $gcd(H_{k,\tau_i}^{(i)}, f_l^{(i)})$ and
$gcd(H_{k,\tau_i}^{(j)}, f_l^{(j)})$ have the same degree.\\

The above very brief outline of the proof already shows that a
detailed proof of this version of the assertion would have involved
a lot of careful, possibly confusing, notation (just imagine a
detailed account of Schoof's Algorithm with this notation!). It
should be clear however, that the proof follows the same string of
arguments.

\section{Advanced topics}

There are some further simplifications of the assertion that I have
made.  Some of the topics underlying the full claim made by Dr
Kayal, and Prof Poonen's addition to this, are somewhat too advanced
to give a ``brief'' explanation of them before being able to prove
the theorem.  Details
of such topics however can be found e.g. in \cite{poonen}, \cite{sil1} and \cite{kobl1}.\\

In its original form, the assertion has the following underlying
 hypothesis:

\begin{quote}
Let $h(t)$ in $\mathbb{Z}[t]$ be a nonzero polynomial, and let $C$
be a smooth projective curve of genus g over $\mathbb{Z}[t] /
(h(t))$.

The hypothesis is: There exists a $C$ as above such that for each
sufficiently large primes $p$, the zeta functions of the curves
$C_{p,\alpha_i}$ are distinct.

\end{quote}
\subsection{Curves of Genus g}

In this dissertation I have restricted myself to the case of very
``simple'' curves, i.e. elliptic curves, which are also often
defined as non-singular curves of genus 1.\\

Giving a detailed discussion about curves of higher genus would take
us too far afield in this dissertation.  As a brief description of
``genus'' however, I will just say that any curve $F(x, y) = 0$ has
a non-negative integer $g$ associated with it; $g$ is referred to as
the genus.  In general (for example if the curve is non-singular),
the genus increases as the degree of $F$ increases.  \\
For curves of the form $Y^2 = F(X)$, we have that $deg(F) = 2g +1$.
Hence in our case the genus is 1, and a curve of the form $Y^2 = X^5
+ a_4X^4 + \dots + a_0$ has genus 2.\\

The addition that Prof Poonen made to Dr Kayal's initial claim is in
fact to
do with such curves of genus greater or equal to 2.\\

In trying to find ``applicable'' curves for this problem, he
remarked that if one used a curve of higher genus, the probability
of the hypothesis being fulfilled would be much higher than for
elliptic curves.\\

In fact, he claimed that the probability of finding an elliptic
curve over a ring $\mathbb{F}_q[z]/(h(z))$, for which the fibres at
each root $\alpha_i$ (that is, the ``subcurves'' $C_i$ for all $i$)
have the same number of
rational points modulo some prime number $l$, is of order $1/q$.\\

For curves of higher genus, that probability is, according to Prof
Poonen, much smaller - in fact it is of order $1/q^2$.  This is of
course an important result if one tries to find curves to
apply this theorem to.\\

We can however find suitable elliptic curves, too, that fulfil our
hypothesis.  For example, consider the curve
\begin{displaymath} C: Y^2 = X^3 + zX \textrm{ over } R =
\mathbb{F}_5[z]/(z^2-1) \, .\end{displaymath}

The subcurves, i.e. $C$ evaluated at the roots of $(z^2 - 1) =
(z+1)(z-1)$ are given by
\begin{displaymath} C_+: Y^2 = X^3 + X \textrm{ and } C_-: Y^2 = X^3
- X \, , \end{displaymath} both defined over $\mathbb{F}_5$.\\

By inspection, we can see that the points in $C_+(\mathbb{F}_5)$ are
\begin{displaymath} \{o, (0,0),(2,0),(3,0)\},\end{displaymath} so there are 4 of them.\\

On the other hand, $C_-(\mathbb{F}_5)$ has the points
\begin{displaymath}\{o, (0,0), (1,0), (2, 1), (2, -1), (3, 2), (3, -2), (4,
0)\}\end{displaymath} - a set of 8
rational points!\\

Clearly, $4 \ne 8$ (mod $3$).  So this curve satisfies our
hypothesis and could be used in applications of the theorem
(although of course it would be a fairly trivial and pointless
example).\\

Now in order to show that what we have proved above is (almost) the
same as the original assertion by Dr Kayal, we will just need to
understand what the ``Zeta function of $C$'' is.

\subsection{The Zeta function} \label{subs:zeta}

Let $C$ be a curve defined over $\mathbb{F}_q$.  Clearly if $C$ is
defined over $\mathbb{F}_q$ then it is also defined over
$\mathbb{F}_{q^n}$ for all $n \ge 1$.  It may therefore be
interesting to consider
\begin{displaymath}
N_n=\#C(\mathbb{F}_{q^n}) \end{displaymath} for $n \ge 1$, i.e. the
number of rational points on $C$ over
$\mathbb{F}_{q^n}$.

\begin{predef}[The Zeta Function] Define the series
\begin{equation}\label{eq:zeta} Z(E;T) = \exp\bigg(\sum_{n \ge 1} \frac{N_n}{n} \cdot
T^n\bigg) \end{equation} for an indeterminate T. This is called the
\textit{Zeta function} of C over $\mathbb{F}_q$.
\end{predef}

Due to work by Hasse - and for a more general case extending to
curves of genus higher than 1, by Weil - we can show that the Zeta
function has a simpler form:

\begin{theorem}[Weil conjectures for an elliptic curve]
Let $C$ be a curve defined over $\mathbb{F}_q$.  Denote by $c_n$ the
trace of Frobenius of $C$ over $\mathbb{F}_{q^n}$, i.e. $c_n =
\#C(\mathbb{F}_{q^n}) - q^n - 1$.  The Zeta function is a rational
function of T and takes the form
\begin{equation}\label{eq:ZetaT}
Z(C;T) = \frac{P(T)}{(1-T)(1-qT)} \end{equation} where $P(T) = 1 -
c_1T + qT^2 = (1 - \alpha)(1 - \bar{\alpha})$.  Furthermore, the
discriminant of $P(T)$ is non-positive and the magnitude of $\alpha$
is $\sqrt{q}$.
\end{theorem}

A proof of this theorem can be found e.g. in \cite{kobl2}, for the
case $g = 1$, or \cite{sil2}, for curves of higher genus. Note in
particular that the last line of the theorem implies that $c_1^2 <
4q$, which is Hasse's Theorem.\\

More importantly however, if we take the derivative of the logarithm
of both sides in (\ref{eq:ZetaT}), substituting in (\ref{eq:zeta})
for the left hand side, we can show after some straightforward
series manipulations and partial fraction expansions that this
implies
\begin{equation*} \label{eq:Nn} N_n = \#C(\mathbb{F}_{q^n}) = q^n + 1
- \alpha^n - \bar{\alpha}^n = |1 - \alpha^n|^2 \, .
\end{equation*}

Since both $c_1$ and $\alpha$ can be immediately derived from
knowledge of $N_1$, we can uniquely determine $N_n$ for all $n \ge
1$ once we know the number of rational points of $C$ over the base
field $\mathbb{F}_q$.\\

It should be clear that this result is an extremely important one
which has many useful applications not only in attempting to prove
the above theorem.  However, if we return to Kayal/Poonen's claim,
we can now simplify the proof of the theorem as follows:\\

The underlying hypothesis for our factorisation is a different Zeta
function for the fibres at $\alpha_i$ (mod $l$) and $\alpha_j$ (mod
$l$) for some prime $l$. Now that we know that we can determine
$N_n$ unambiguously from computing $N_1$, and since $Z(C;T)$ only
depends on the $N_n$, it will be sufficient to use a curve for which
the fibres at $\alpha_i$ and $\alpha_j$ have a different number of
rational points over $\mathbb{F}_q$, modulo $l$.  This is clearly
what we have done above.

\chapter{Running Time} \label{ch:running}

To finish this dissertation off, I will now give an account of the
running time of the
algorithm.\\

Since algorithms and computations thereof is a rather broad
mathematical subject on its own, and one which I assume the reader
to be unfamiliar with, I will restrict this discussion to the key
points.  I recommend in particular \cite{bach} to the interested
reader; it contains a comprehensive introduction to this subject,
and will also fill in some details about the running time of the
individual steps in our
algorithm that I will omit.\\

\section{Introduction to running times}

In order to calculate the running time of an algorithm, we count the
number of basic operations performed by the algorithm on the
``worst-case input''.  The \textit{worst-case input} is the input
for which the most basic operations are required.  We count the
basic operations as follows:

\subsection{Definition}

Let $n \in \mathbb{Z}$.  Define
\begin{equation*} lg \, n:= \left\{
\begin{array}{ll}
1, & \textrm{if } n = 0; \\
1 + \lfloor log_2 |n| \rfloor , & \textrm{if } n \ne 0. \end{array}
\right.
\end{equation*}

Then $lg \, n$ counts the number of bits in the binary
representation of $n$.\\

A \textit{step} is the fundamental unit of computation.  Now,
different situations require different units.  For example,
analysing a sorting algorithm would require counting the number of
comparison steps, whereas in the case of an algorithm that computes
the evaluation of a polynomial at a certain point we may want to
count each addition, subtraction and multiplication as a
single step.\\

In general we therefore adapt the convention to equate ``step'' with
``bit operation'':  We write all integers in binary code, so we are
only working with variables that take the values 0 or 1.  We then
perform logical operations on these variables:  conjunction
($\land$), disjunction ($\lor$) and negation ($\sim$).  Each of those operations takes 1 bit.\\

The \textit{running time}, or \textit{cost of computation}, is then
the total number of such logical operations performed in an
algorithm. It
depends on the size of the input.\\

For example, the operation $a + b$ takes $lg \, a + lg \, b$ bit
operations.  However we usually just aim to find an upper bound of
the running time, rather than an exact number; we therefore only
note that the running time is $O(lg \, a + lg \, b)$ (where the $O$
is the ``Big-Oh-notation'', which should be well-known to the
reader).

\begin{Lemma} \label{runningtime}

If we have a sequence of operations in an algorithm, say $P$ and
$Q$, then we have that the running time of the algorithm
``operation $P$ followed by operation $Q$'' is\\

$Time(P \, ; Q) = Time(P) + Time(Q)$.\end{Lemma}

\section{Running time of our algorithm}

Lemma \ref{runningtime} tells us that in order to compute the
precise running time of the algorithm, we need to add up the running
times required for each individual part of the algorithm.\\

However, to find an upper bound of the running times, it will be
sufficient to find an upper bound for the part of the algorithm that
has the largest running time, as stated in the following Lemma:

\begin{Lemma} \label{bigolemma}
If $f(x) = O(g(x))$ then $f(x) + g(x) = O(g(x))$.\end{Lemma}

Let us now go through the steps of Schoof's Algorithm, and analyse
the amount of computations involved in each step.  I will use
several standard results for running times without proof;
details can be found e.g. in \cite{bach}.\\

It is a well-known result that there exists a universal constant $C$
such that
\begin{displaymath} \prod_{\substack{l \textrm{ prime} \\ 2\le l \le l_{max}}} l > C \cdot
e^L.
\end{displaymath}

for every $L>0$. For a proof of this, see e.g. \cite{ros}.\\

So we can take $l_{max}:= O(lg \, q)$. The number of primes
occurring in the product is $O(lg \, q)$ and the primes $l$
themselves are clearly also $O(lg \, q)$.\\

Now consider the running time involved in Step (2) of the algorithm.
This step clearly requires the largest amount of computation, so
that its
running time will in the end dominate over the others.\\

We now need to state some more results from complexity theory.

\begin{predef}

Let $f(x) \in R[X]$ for some ring $R$ with $|R| = p^m$.  Define
\begin{equation*} lg \, f:= \left\{
\begin{array}{ll}
1, & \textrm{if } f = 0; \\
(1 + deg (f)) lg \, |R|, & \textrm{if } f \ne 0.
\end{array} \right.
\end{equation*}\end{predef}

\begin{Lemma}

Let $f$, $g$ be polynomials in $R[X]$. Then

\begin{enumerate}
\item  $f \pm g$ can be
computed with $O(lg \,f + lg \, g)$ bit operations.
\item $f \cdot g$ can be computed using $O((lg \, f)(lg \, g))$ bit
operations.
\item Computing the greatest common divisor of $f$ and $g$ also
requires $O((lg \, f)(lg \, g))$ bit operations.
\end{enumerate}
\end{Lemma}

From Corollary \ref{fdef}, Part \ref{degf}, we have that $deg(f_l) =
O(l^2)$,
and from above we know that $l_{max} = O(lg \, q)$.  \\

Computing the $H_{k, \tau}$ will involve computing $x^q, y^q,
x^{q^2}$ and $y^{q^2}$ (reduced modulo the curve equation) modulo
$f_l$.\\

For $x^q$ and $x^{q^2}$, this will require $O(lg \, q)$
multiplications in the ring each - hence $O((lg \, q)^2)$ together.
Reducing modulo $f_l$ takes $O(l^2) = O((lg \, q)^2)$ bit
operations, so the computation of $x^q$ and $x^{q^2}$ will require
$O((lg \, q)^4)$ multiplications in the ring.\\

Since the order of the ring in our case is $|\mathbb{F}_q| = O(lg \,
q)$, multiplication of any two elements in $R$ takes
$O((lg \, q)^2)$ bit operations.\\

We therefore need $O((lg \, q)^6)$ bit operations in total to
compute $x^q$ and $x^{q^2}$.  For $y^q$ and $y^{q^2}$, the
computations are similar and hence their complexity
will not affect the asymptotic upper bound.\\

Now the $x^q, y^q, x^{q^2}$ and $y^{q^2}$ are computed once for each
prime $l$, so $O(lg \, q)$ times, and then stay the same for each
$\tau$.  Now $\tau$ is also $O(lg \,q)$, so we have $O((lg \, q)^7)$
bit operations for each prime $l$. \\

Finally, we have $O(lg \, q)$ primes $l$, so the complexity in the
entire Step (2) of Schoof's
Algorithm amounts to $O((lg \, q)^8)$ bit operations.\\

This does indeed dominate the computations of both $l_{max}$ and the
Chinese Remainder Theorem in the last step, so we will not have to
compute the complexities involved in those (we are not concerned
with the latter anyway though, as in \textit{our} algorithm, we will
never get as far as computing the group order!).\\

In fact, one can make improvements to find a slightly lower ``upper
bound'' for the complexity, but let us finish this dissertation with
the conclusion that our algorithm computes the factors of $h(z)$
over $\mathbb{F}_q$ using at most $O((lg \, q)^8)$ bit operations.\\

If we convert our ``Big-Oh-notation'' to a polynomial, we can say
that, indeed, the running time of the algorithm is bounded by a
polynomial in
$log \, q$, as asserted at the beginning of this dissertation.\\

\chapter{Acknowledgements}

I would like to thank Dr Lauder, both for drawing my attention to
the idea for this project, and for all his support throughout the
last two terms.

I would also like to express my gratitude to Dr Flynn, who in his
lecture course on Elliptic Curves has provided me with a strong
background in the study of this field and
whose ideas for my dissertation have been extremely helpful.\\

The insightful email conversations with Prof Poonen have been
invaluable, and I
am very grateful for all his inspirations.\\

The online guide to \LaTeX  provided by the Mathematical Institute
(written by Tobias Oetiker) has been a fantastic way for me to teach
myself in a very limited amount of time how to use this language.
Also the resources provided by the University Library Services and
the Computing Services ought to be mentioned here - in terms of
availability of books, software and general help, they have been
extremely efficient.\\

Last but certainly not least, I am extremely grateful to my friends
for their ``moral support''.  \\

In particular, I thank my parents for all their advice and care.
Without it, the final few weeks would have been unthinkable.

\begin{center} \textit{Thank you.}\\

L.D.S.
\end{center}

\end{document}